\documentclass[a4paper,12pt]{amsart}
\usepackage{amssymb}
\usepackage{ifthen}
\usepackage{graphicx}
\usepackage{color}
\nonstopmode
\newtheorem{thm}{Theorem}
\newtheorem{rem}{Remark}
\newtheorem{cor}{Corollary}
\newtheorem{lem}{Lemma}
\newtheorem{prop}{Proposition}

\newtheorem{claim}{Claim}

\newtheorem{conj}{Conjecture}
\newtheorem{prob}{Problem}

\theoremstyle{definition}
\newtheorem{defn}{Definition}
\newtheorem{example}{Example}

\newtheorem{case}{Case}

\newenvironment{pf}[1][]{%
 \vskip 1mm
 \noindent
 \ifthenelse{\equal{#1}{}}%
  {{\slshape Proof. }}%
  {{\slshape #1.} }%
 }%
{\qed\medskip}

\newcounter{alphabet}
\newcounter{tmp}
\newenvironment{Thm}[1][]{\refstepcounter{alphabet}%
\bigskip%
\noindent%
{\bf Theorem \Alph{alphabet}}%
\ifthenelse{\equal{#1}{}}{}{ (#1)}%
{\bf .} \itshape}{\vskip 8pt}

\makeatletter
\newcommand{\Ref}[1]{\@ifundefined{r@#1}{}{\setcounter{tmp}{\ref{#1}}\Alph{tmp}}}
\makeatother

\newenvironment{Lem}[1][]{\refstepcounter{alphabet}%
\bigskip%
\noindent%
{\bf Lemma \Alph{alphabet}}%
{\bf .} \itshape}{\vskip 8pt}

\newcommand{\IR}{{\mathbb R}}

\newcommand{\IC}{{\mathbb C}}
\newcommand{\ID}{{\mathbb D}}

\newcommand{\real}{{\operatorname{Re}\,}}

\def\be{\begin{equation}}
\def\ee{\end{equation}}

\newcommand{\bee}{\begin{enumerate}}
\newcommand{\eee}{\end{enumerate}}

\newcommand{\blem}{\begin{lem}}
\newcommand{\elem}{\end{lem}}
\newcommand{\bthm}{\begin{thm}}
\newcommand{\ethm}{\end{thm}}
\newcommand{\bcor}{\begin{cor}}
\newcommand{\ecor}{\end{cor}}
\newcommand{\beg}{\begin{example}}
\newcommand{\eeg}{\end{example}}
\newcommand{\bprop}{\begin{prop}}
\newcommand{\eprop}{\end{prop}}
\newcommand{\begs}{\begin{examples}}
\newcommand{\eegs}{\end{examples}}
\newcommand{\bdefe}{\begin{defn}}
\newcommand{\edefe}{\end{defn}}
\newcommand{\bprob}{\begin{prob}}
\newcommand{\eprob}{\end{prob}}
\newcommand{\bques}{\begin{ques}}
\newcommand{\eques}{\end{ques}}
\newcommand{\bei}{\begin{itemize}}
\newcommand{\eei}{\end{itemize}}

\newcommand{\bde}{\begin{deter}}
\newcommand{\ede}{\end{deter}}
\newcommand{\bca}{\begin{case}}
\newcommand{\eca}{\end{case}}
\newcommand{\bcl}{\begin{claim}}
\newcommand{\ecl}{\end{claim}}
\newcommand{\bcon}{\begin{conj}}
\newcommand{\econ}{\end{conj}}
\newcommand{\bcons}{\begin{conjs}}
\newcommand{\econs}{\end{conjs}}
\newcommand{\br}{\begin{rem}}
\newcommand{\er}{\end{rem}}
\newcommand{\brs}{\begin{rems}}
\newcommand{\ers}{\end{rems}}
\newcommand{\bo}{\begin{obser}}
\newcommand{\eo}{\end{obser}}
\newcommand{\bos}{\begin{obsers}}
\newcommand{\eos}{\end{obsers}}
\newcommand{\bpf}{\begin{pf}}
\newcommand{\epf}{\end{pf}}
\newcommand{\ba}{\begin{array}}
\newcommand{\ea}{\end{array}}
\newcommand{\beq}{\begin{eqnarray}}
\newcommand{\beqq}{\begin{eqnarray*}}
\newcommand{\eeq}{\end{eqnarray}}
\newcommand{\eeqq}{\end{eqnarray*}}

\date{May 2018}

\begin{document}
\title[Rotations and convolutions of  harmonic convex mappings]
{Rotations and convolutions of harmonic convex mappings}


\author{Liulan Li and Saminathan  Ponnusamy
}

\address{Liulan Li, College of Mathematics and Statistics
 (Hunan Provincial Key Laboratory of Intelligent Information Processing and Application),
Hengyang Normal University, Hengyang,  Hunan 421002, People's
Republic of China} \email{lanlimail2012@sina.cn}

\address{S. Ponnusamy,
 Department of Mathematics,
Indian Institute of Technology Madras, Chennai-600 036, India.}
\email{samy@iitm.ac.in}



\maketitle

%
%

\begin{abstract}
In this paper, we consider the convolutions of slanted half-plane
mappings and strip mappings and generalize related results in
general settings. We also consider a class of harmonic mappings
containing slanted half-plane mappings and strip mappings and, as a consequence, we prove that
the convex combination of such mappings is convex.
\end{abstract}

\maketitle \pagestyle{myheadings}
\markboth{Liulan Li and Saminathan Ponnusamy}{Rotations and convolutions of harmonic mappings}

\section{Introduction and Preliminary results}


In this article, we will consider complex-valued harmonic mappings
$f$ defined on the open unit disk ${\mathbb D}=\{z \in {\mathbb
C}:\, |z|<1\}$, which have the canonical representation of the form
$f=h+\overline{g}$, where $h$ and $g$ are analytic in $\ID$. This
representation is unique with the condition $g(0)=0$. In terms of
the canonical decomposition of $f$, the Jacobian $J_f$ of
$f=h+\overline{g}$ is given by $J_f(z)=|h'(z)|^2 -|g'(z)|^2$.
According to Inverse Mapping Theorem, if the Jacobian of a $C^1$
mapping from $\ID$ to $\IC$ is different from zero, then the
function is locally univalent in $\ID$. The classical result of Lewy implies
that the converse of this statement also holds for harmonic
mappings.  Thus, every harmonic mapping $f$ on $\ID$ is locally
one-to-one and sense-preserving on $\ID$ if and only if $J_f(z)>0$
in $\ID$, i.e. $|h'(z)|>|g'(z)|$ in $\ID$. For locally
one-to-one mappings, the condition $J_f(z)> 0$
is equivalent to the existence of an analytic function $\omega _f $
in $\ID$ such that
$$ |\omega _f (z)|<1 ~\mbox{ for }~z\in \ID,
$$
where $\omega _f(z)=g'(z)/h'(z)$ is called the dilatation of $f$.
When there is no confusion,  it is often convenient to use $\omega $
instead $\omega _f$. Let ${\mathcal H}=\{f=h+\overline{g}:\,
~\mbox{$h(0)=g(0)=0$ and $h'(0)=1$} \}$. The class ${\mathcal H}_0$
consists of those functions $f\in {\mathcal H}$ with $g'(0)=0$.

The family of all sense-preserving univalent harmonic mappings from
${\mathcal H}$ will be denoted by ${\mathcal S}_H$, and let
$\mathcal{S}_{H}^{0}={\mathcal S}_H\cap {\mathcal H}_0$. Clearly,
the familiar class ${\mathcal S}$ of normalized analytic univalent
functions in $\ID$ is contained in ${\mathcal S}_H^{0}$. The class
${\mathcal S}_H$ together with its geometric subclasses has been
studied extensively  by Clunie and Sheil-Small
\cite{Clunie-Small-84} and investigated subsequently by several
others (see \cite{Du} and the survey article \cite{SaRa2013}). In
particular, we consider the convolution properties of the class
${\mathcal K}_H$ (resp.  ${\mathcal K}_H^0$) of functions ${\mathcal
S}_H$ (resp. ${\mathcal S}_H^0$) that map the unit disk $\ID$ onto a
convex domain. The class  ${\mathcal K}_H^0$ has been studied
extensively. See \cite{Clunie-Small-84,Du,Good02,SaRa2013} and also
the recent articles \cite{FHM,HM,LiuPo}.


The main aim of this article is to derive convolution results and in particular, we extend
and revaluate many of the recent results on slanted half-plane and strip mappings.
The main  results are Theorems \ref{LiS16-th1}, \ref{LiS16-th2a} and \ref{unit1}.
In order to achieve these results, we need some preparation. We begin to include several preliminary results, observations
and several useful lemmas  in the following subsections.

\subsection{Slanted half-plane and strip mappings}
We say that a domain $D$ is called convex in the direction $\alpha$ $(0\leq \alpha<
\pi)$ if the intersection of $D$ with each line parallel to
the line through $0$ and $e^{i\alpha}$ is an interval or the empty
set. The class of univalent harmonic functions $f$ for which the range $D=f(\ID)$ is a
domain that is convex in the direction $\alpha$
plays a crucial role in deriving properties of convolution of harmonic mappings.
Such functions are called convex in the direction $\alpha$. See, for example,
\cite{Clunie-Small-84,Do,DN,HengSch70,HengSch73,M,RZ}. Functions that are convex in
the direction $\alpha =0$ (resp. $\alpha =\pi/2$) is referred to as
convex in the real (resp. vertical) direction.

It is known that \cite{Clunie-Small-84}, a harmonic mapping
$f=h+\overline{g}$ belongs to $\mathcal{K}_{H}^{0} $ if and only if, for each $\alpha \in [0,
\pi )$, the analytic function $F=h-e^{2i\alpha}g$ belongs to $\mathcal{S}$
and is convex in the direction $\alpha$. This result is instrumental in deriving many properties of
the class $\mathcal{K}_{H}^{0}$ by transferring information from conformal case.

\bdefe\label{def1} A function  $f=h+\overline{g} \in {\mathcal S}_H$
is said to be a slanted half-plane mapping with $\gamma$
($0\leq\gamma<2\pi$) and $a\in \ID$ if $f$ maps $\ID$ onto the
half-plane $H_\gamma^a :=\left \{w:\,{\rm Re}\left
(\frac{e^{i\gamma}}{1+a}w\right ) >-\frac{1}{2}\right \}$. The class
of all such mappings is denoted by  ${\mathcal S}(H^a_{\gamma})$.
\edefe

In \cite{LiPo3,LiuPo}, this definition was considered for $a\in (-1,1)$.
In view of the Riemann mapping theorem, it is now a standard procedure to obtain the following.

\bprop \label{LP16Prop1}
Each $f=h+\overline{g} \in{\mathcal S}(H^a_{\gamma})$ can be written as
\begin{equation}\label{li4-eq2}
h(z)+e^{-2i(\gamma +\gamma _a)}g(z)=\frac{(1+a')z}{1-e^{i(\gamma +\gamma _a)}z},
\end{equation}
where  $\gamma _a=\arg (1+\overline{a})$ and $g'(0)=a'e^{2i(\gamma +\gamma _a)}$ with $a'\in \ID$ and $a'=|1+a|-1$.
\eprop
\bpf
In fact, \eqref{li4-eq2} follows from writing
$$\frac{e^{i\gamma}}{1+a}f=\frac{e^{i(\gamma +\gamma _a)}}{|1+a|}(h+\overline{g})=
\frac{e^{i(\gamma +\gamma _a)}h +\overline{e^{-i(\gamma +\gamma
_a)}g}}{|1+a|}
$$
so that $f=h+\overline{g} \in{\mathcal S}(H^a_{\gamma})$ is equivalent to
$${\rm Re}\left (\frac{e^{i\gamma}}{1+a}f(z)\right )={\rm Re}\left (\frac{e^{i(\gamma +\gamma _a)}(h(z) +e^{-2i(\gamma +\gamma _a)}g(z))}{|1+a|}\right )>-\frac{1}{2}, \quad z\in\ID.
$$
Note that $h(0)=g(0)=h'(0)-1=0$ and
$$g'(0)=(|1+a|-1)e^{2i(\gamma +\gamma _a)}=\left((1+\overline{a})e^{-i\gamma _a} -1\right)e^{2i(\gamma +\gamma _a)}=
\left((1+\overline{a})  -e^{i\gamma _a}\right)e^{i(2\gamma +\gamma
_a)}.
$$
Also, we remark that $g'(0)\in \ID$, because $|g'(0)|<1$ is
equivalent to the inequality $0<|1+a|<2$ which is clearly true. The
Riemann mapping theorem gives the desired representation
\eqref{li4-eq2}. \epf

 Obviously, each $f\in {\mathcal S}(H^a_{\gamma})$ (resp.
$\mathcal{S}(H^0_{\gamma})$) belongs to the convex family
$\mathcal{K}_{H}$ (resp. $\mathcal{K}_{H}^{0}$). We remark that
there are infinitely many \textit{slanted half-plane mappings} with
a fixed $\gamma$ and a fixed $a$.

For notational consistency, it is appropriate to consider harmonic mappings
$f=h +\overline{g}\in {\mathcal S}(H_\gamma^a)$ with the dilatation
$$\omega (z)=-e^{2i(\gamma +\gamma _a)}\frac{e^{i(\gamma +\gamma _a)}z-a'}{1-a'e^{i(\gamma +\gamma _a)}z},
~\mbox{ and
}~h(z)+e^{-2i(\gamma +\gamma _a)}g(z)=\frac{(1+a')z}{1-e^{i(\gamma +\gamma _a)}z}.
$$
Solving these two equations with the help of the fact that $\omega
=g'/h'$ gives slanted half-plane mappings with $\gamma$ and $a$ as
\be\label{LiS16-eq12} f^a_\gamma=h^a_\gamma+\overline{g^a_\gamma},
\ee where
$$h^a_\gamma(z)=\frac{z-\frac{1+a'}{2}e^{i(\gamma +\gamma _a)}z^2}{(1-e^{i(\gamma +\gamma _a)}z)^2}=\frac{(1+a')I(z)+(1-a')zI'(z)}{2}, \quad
I(z)=\frac{z}{1-e^{i(\gamma +\gamma _a)}z},
$$
and
$$e^{-2i(\gamma +\gamma _a)}g^a_\gamma(z)=\frac{a'z-\frac{1+a'}{2}e^{i(\gamma +\gamma _a)}z^2}{(1-e^{i(\gamma +\gamma _a)}z)^2}=
\frac{(1+a')I(z)-(1-a')zI'(z)}{2} .
$$
Obviously, when $a$ is real and $a\in (-1,1)$, we have $\gamma_a=0$ and $a'=a$ so that $f^a_\gamma$ coincides with the one
considered in \cite{LiPo3}.

Moreover, if $a\in(-1,1)$, then  $f=h+\overline{g} \in{\mathcal S}(H^a_{\gamma})$ satisfies the relation
$$h(z)+e^{-2i\gamma}g(z)=\frac{(1+a)z}{1-e^{i\gamma}z}, \quad z\in\ID.
$$
It is worth recalling that functions $f\in{\mathcal S}(H^a_{\gamma})$ with $\gamma=0$ and $a\in
(-1,1)$ are usually referred to as \textit{the right half-plane
mappings}, especially when $a=0$. In the later case, we have  $f_{0}=h_{0}+\overline{g_{0}}\in {\mathcal S}(H^0_{0})$
with the dilatation $\omega_{0}(z)=-z$,  where
\be\label{li4-eq2b}
 h_{0}(z)=\frac{z-\frac{1}{2}z^2}{(1-z)^2}=
\frac{1}{2}\left (\frac{z}{1-z}+\frac{z}{(1-z)^2} \right )
\ee
and
\be\label{li4-eq2c}
g_{0}(z)=\frac{-\frac{1}{2}z^2}{(1-z)^2}=\frac{1}{2}\left
(\frac{z}{1-z}-\frac{z}{(1-z)^2} \right ),
\ee
so that
$$
h_{0}(z)+g_{0}(z)=\frac{z}{1-z}.
$$
The function $f_0$ is central for many extremal problems for the convex family $\mathcal{K}_{H}^{0}$, and
$f_0$ maps the unit disk onto the half-plane $\{w:\,{\rm Re} w  >-1/2\}$.

\bdefe
For $0<\beta<\pi$ and $b\in \ID$, let ${\mathcal S}(\Omega^b_\beta)$ denote the class of functions $f$
from ${\mathcal S}_H$ such that $f$ maps $\ID$ onto the strip
$\Omega^b_\beta,$  where
$$\Omega^b_\beta :=\left \{w:\, \frac{\beta-\pi}{2\sin\beta}<
 {\rm Re}\left (\frac{1}{1+b}w\right )<\frac{\beta}{2\sin\beta}\right
\}.
$$
\edefe

As with Definition \ref{def1} and Proposition \ref{LP16Prop1}, it is a simple exercise to obtain the
following:

\bprop \label{LP16Prop1a}
Each $f=h+\overline{g} \in {\mathcal S}(\Omega^b_\beta)$ has the form
\be\label{li4-eq3}
h(z)+e^{-2i\gamma _b}g(z)=\psi (z),
\quad \psi(z)= \frac{(1+b')e^{-i\gamma _b}}{2i\sin\beta}\log\left(\frac{1+ze^{i(\beta +\gamma
_b)}}{1+ze^{-i(\beta-\gamma _b)}} \right),
\ee
where $\gamma _b=\arg (1+\overline{b})$ and $b'=|1+b|-1$.
\eprop

Note that $h(0)=g(0)=h'(0)-1=0$ and
$$g'(0)=(|1+b|-1)e^{2i\gamma _b}=\left((1+\overline{b})e^{-i\gamma _b} -1\right)e^{2i\gamma _b}=
\left((1+\overline{b})  -e^{i\gamma _b}\right)e^{i\gamma _b}.
$$
It follows that each $f\in {\mathcal S}(\Omega^b_\beta)$ (resp.
${\mathcal S}(\Omega^0_\beta)$) belongs to the convex family
$\mathcal{K}_{H}$ (resp. $\mathcal{K}_{H}^{0}$). If $b\in(-1,1)$,
then $b'=b$ so that the last relation for $f=h+\overline{g} \in
{\mathcal S}(\Omega^b_\beta)$ takes the form \be\label{li4-eq3c}
h(z)+g(z)=\frac{1+b}{2i\sin\beta}\log\left(\frac{1+ze^{i\beta}}{1+ze^{-i\beta}}
\right), \ee where $g'(0)=b$. The class of the strip mappings with $b=0$ has
been considered extensively.

\subsection{Convolution of harmonic mappings}
If $f=h+\overline{g}$ and $F=H+\overline{G}$ are two harmonic mappings of the unit disk $\ID$ with power series of the form
$$f(z)
=\sum_{n=1}^{\infty}a_{n}z^{n}+\sum_{n=1}^{\infty}\overline{b_{n}}\overline{z}^{n}
~\mbox{ and }~
F(z)
=\sum_{n=1}^{\infty}A_{n}z^{n}+\sum_{n=1}^{\infty}\overline{B_{n}}\overline{z}^{n},
$$
then the harmonic convolution (or Hadamard product) of $f$ and $F$, denoted by $f*F$, is defined by
$$f*F=h*H+\overline{g*G}=\sum_{n=1}^{\infty}a_{n}A_{n}z^{n}+\sum_{n=1}^{\infty}\overline{b_{n}B_{n}}\overline{z}^{n}.
$$
In particular, the space ${\mathcal H}$ is closed under the operation
$\ast$, i.e. ${\mathcal H} \ast {\mathcal H}\subset {\mathcal H}$,
and $f*F=F*f$. Moreover,
$$(h*H)'(z)=\frac{h(z)}{z}*H'(z), ~\mbox{ i.e. }~z(h*H)'(z)=h(z)*zH'(z)=zh'(z)*H(z),
$$
which will be used in the proof of our theorems.

In the case of conformal mappings, many important properties of
convolution are established (see \cite{samy95,PoSi96,rs1}). For example,
convolution of two univalent convex (analytic) functions is convex.
On the other hand, most of such results do not carry over to the
case of univalent harmonic mappings in $\ID$ (see
\cite{Do,DN,Good02,LiPo1,LiPo2,LiPo3}) even in simplest cases.
For instance, convolution of two mappings from $\mathcal{K}_{H}^{0}$
is not necessarily locally univalent even in restricted cases. In
spite of such drawbacks, some properties of the
convolution of harmonic univalent mappings have been achieved and
these results help to understand how difficult it is to generalize the convolution
properties of harmonic mappings. See \cite{Do,DN,Good02,LiPo1,LiPo2}.

Our main results  are related to the theory of harmonic convolutions which extend
and revaluate many of the recent results.

\subsection{Rotations}
Motivated by the work in \cite{FHM}, for $f=h+\overline{g}\in{\mathcal H}$ and $\mu\in\partial\ID$,  consider the standard
rotation $f^\mu$ of $f$ by the relation
\be\label{LiS16-eq1a}
f^\mu(z)=\overline{\mu}f(\mu z)
\ee
and thus, has the canonical decomposition
\be\label{LiS16-eq1}
f^\mu =h^\mu+\overline{g^\mu }, ~\mbox{ where }~
h^\mu(z)= \overline{\mu}h(\mu z) ~\mbox{ and }~ g^\mu(z)= \mu g(\mu z).
\ee
Since $(h^\mu)'(z)=h'(\mu z)$ and $(g^\mu)'(z)= \mu^2 g'(\mu z)$,  it follows easily that
 \be\label{LiS16-eq2}
 \omega_{f^\mu}(z)=\mu^2\frac{g'(\mu z)}{h'(\mu z)}=\mu^2 \omega_{f}(\mu z).
\ee
Thus, we may formulate the above discussion as
\begin{lem}\label{the dilatation after rotation}
Let $\mu\in\partial\ID$ and $f=h+\overline{g}\in{\mathcal H}$. Then
$f^\mu$ is locally one-to-one and sense-preserving on $\ID$ if and only if $f$ is locally one-to-one
and sense-preserving on $\ID$.
\end{lem}

Moreover, we can easily obtain the following two lemmas.

\begin{lem}\label{the convolution after rotation}
Let $f_j\in{\mathcal H}$, and  $\mu_j\in\partial\ID$, for $j=1,2$.
Then
$$(f^{\mu_1}_1*f^{\mu_2}_2)(z)=\overline{\mu_1 \mu_2}\,(f_1*f_2) (\mu_1 \mu_2 z)=(f_1*f_2)^{\mu_1 \mu_2}(z).
$$
\end{lem}
\bpf Follows from the definitions of $f^{\mu}$ and the convolution.
\epf

The relations \eqref{LiS16-eq1} and \eqref{LiS16-eq2} and the simple fact that
$$\frac{e^{i\gamma}}{1+a}f( e^{i\xi}z) = \frac{e^{i\gamma}e^{i\xi}}{1+a}[e^{-i\xi}h( e^{i\xi}z) +\overline{e^{i\xi}g( e^{i\xi}z) }]
=\frac{e^{i(\gamma +\xi)}}{1+a}f^{ e^{i\xi}}(z)
$$
give the following.

\begin{lem}\label{the slanted half-plane mappings after rotation}
Let $a\in \ID$, $a'=|1+a|-1$,   and $f=h+\overline{g}\in {\mathcal
S}(H^a_{\gamma})$. Then $f^{e^{-i(\gamma+\gamma_a)}}\in {\mathcal
S}(H^{a'}_{0})$, where $\gamma _a=\arg (1+\overline{a})$. Moreover,
if $f=h+\overline{g}\in {\mathcal S}(H^a_{\gamma})$ with
\be\label{LiS16-eq3} h(z)+e^{-2i(\gamma +\gamma
_a)}g(z)=\frac{(1+a')z}{1-e^{i(\gamma +\gamma _a)}z} ~\mbox{ and }~
\omega_{f}
(z)=e^{2i(\gamma+\gamma_a)}\frac{a'-ze^{i(\gamma+\gamma_a)}}{1-a'ze^{i(\gamma+\gamma_a)}},
\ee then $f^{e^{-i(\gamma+\gamma_a)}}=H+\overline{G}\in {\mathcal
S}(H^{a'}_{0})$ with
\be\label{LiS16-eq4} H(z)+
G(z)=\frac{(1+a')z}{1-z}~\mbox{ and }~
\omega_{f^{e^{-i(\gamma+\gamma_a)}}} (z)=\frac{a'-z}{1-a'z}. \ee

If $f=h+\overline{g}\in {\mathcal S}(H^a_{\gamma})$ with
$$\omega_{f}
(z)=e^{2i(\gamma+\gamma_a)}\frac{a'+ze^{i\theta}}{1+a'ze^{i\theta}},
$$
then $f^{e^{-i(\gamma+\gamma_a)}}=H+\overline{G}\in {\mathcal
S}(H^{a'}_{0})$ with
$$\omega_{f^{e^{-i(\gamma+\gamma_a)}}}
(z)=\frac{a'+ze^{i(\theta-\gamma-\gamma_a)}}{1+a'ze^{i(\theta-\gamma-\gamma_a)}}.
$$
\end{lem}

\begin{lem}\label{the strip mappings after rotation}
Let $0<\beta<\pi$, $b\in \ID$, $b'=|1+b|-1$ and
$f=h+\overline{g}\in{\mathcal S}(\Omega^b_\beta)$. Then
$f^{e^{-i\gamma_b}} \in{\mathcal S}(\Omega^{b'}_\beta)$, where
$\gamma _b=\arg (1+\overline{b})$. Moreover, if
$f=h+\overline{g}\in{\mathcal S}(\Omega^b_\beta)$ with
$$h(z)+e^{-2i\gamma
_b}g(z)=\frac{(1+b')e^{-i\gamma
_b}}{2i\sin\beta}\log\left(\frac{1+ze^{i(\beta +\gamma
_b)}}{1+ze^{-i(\beta-\gamma _b)}} \right) ~\mbox{ and }~ \omega_{f}
(z)=e^{2i\gamma_b}\frac{b'+ze^{i\theta}}{1+b'ze^{i\theta}},
$$
then $f^{e^{-i\gamma_b}}=H+\overline{G}\in {\mathcal
S}(\Omega^{b'}_\beta)$ with
$$ H(z)+
G(z)=\frac{(1+b')}{2i\sin\beta}\log\left(\frac{1+ze^{i\beta
}}{1+ze^{-i\beta}} \right)~\mbox{ and }~ \omega_{f^{e^{-i\gamma_b}}}
(z)=\frac{b'+ze^{i(\theta-\gamma _b)}}{1+b'ze^{i(\theta-\gamma
_b)}}.
$$
\end{lem}

\subsection{Basic results on harmonic convolution}

In 2001,  Dorff in \cite{Do} discussed the convolution of functions from
${\mathcal S}(H^0_{0})$ with that of functions from ${\mathcal S}(\Omega^0_\beta)$ in the following form.

\begin{Thm}{\rm (\cite[Theorem 7]{Do})}\label{ThmB} If $f_1\in
{\mathcal S}(H^0_{0})$, $f_2\in{\mathcal S}(\Omega^0_\beta)$ and
$f_1\ast f_2$ is locally univalent in $\ID$, then $f_1\ast f_2$ is
convex in the real direction.
\end{Thm}

In 2012,  Dorff et al. in \cite{DN} proved the following result which concerns the convolution of a function
from ${\mathcal S}(H^0_{\gamma_1})$ with that of function from ${\mathcal S}(H^0_{\gamma_2})$.

\begin{Thm}{\rm (\cite[Theorem 2]{DN})}\label{ThmA}
If $f_k\in {\mathcal S}(H^0_{\gamma_k})$, $k=1, 2$, and $f_1\ast
f_2$ is locally univalent in $\ID$, then $f_1\ast f_2$ is convex in
the direction $-(\gamma_1 +\gamma_2)$.
\end{Thm}

By a similar reasoning as in \cite{DN}, we generalize Theorem
\Ref{ThmA} to the setting ${\mathcal S}(H^a_{\gamma})$ for
$a\in(-1,1)$. More precisely, we have

\begin{Thm}\label{convolution of slanted half-plane}
{\rm (\cite{LiPo3})}
Let $a_k\in(-1,1)$ and  $f_k\in {\mathcal
S}(H^{a_k}_{\gamma_k})$ for $k=1,2$. If $f_1\ast f_2$ is locally
univalent in $\ID$, then $f_1\ast f_2$ is convex in the direction
$-(\gamma_1 +\gamma_2)$.
\end{Thm}

In the next section, we will give a natural generalization of Theorem
\Ref{ThmB}. In order to prove it, we need the following theorems.

\begin{Thm}{\rm (\cite[Lemma 2.7]{rs1})}\label{ThmC}
Let $\phi(z)$ and $\Psi(z)$ be analytic in $\ID$ with $\phi(0) = \Psi(0)= 0$.
If $\phi(z)$ is convex and $\Psi(z)$ is starlike, then for each
function $F(z)$ analytic in $\ID$ and satisfying $\real F(z) > 0$,
we have
$$\real\left\{\frac{(\phi*F\Psi)(z)}{(\phi*\Psi)(z)}\right\}>0,\;\ \forall z\in\ID.
$$
\end{Thm}

Here we say that an analytic function $\Psi$ such that $ \Psi(0)= 0$  is starlike in $\ID$ if  $\Psi(z)$ maps $\ID$ univalently onto a domain which is
starlike with respect to the origin, i.e., $tw\in \Psi(\ID)$ whenever $w\in \Psi(\ID)$ and $t\in [0,1]$.

\begin{Thm}\label{ThmD}
{\rm (\cite{Clunie-Small-84})}
A harmonic $f=h+\overline{g}$ locally univalent in $\ID$ is a
univalent mapping of $\ID$ onto a domain convex in the direction
$\alpha$ $(0\leq \alpha <\pi)$ if and only if $h-e^{2i\alpha}g$ is a
conformal univalent mapping of $\mathbb{D}$ onto a domain convex in
the direction $\alpha$.
\end{Thm}

\begin{Thm}\label{ThmE}{\rm (\cite[Theroem 1]{RZ})}
Let $\varphi(z)$ be a non-constant function regular in $\ID$. The
function $\varphi(z)$ maps $\ID$ univalently onto a domain convex in
the real direction if and only if there are numbers $\mu$ and $\nu$,
$0\leq \mu <2\pi$ and $0\leq \nu\leq \pi$, such that
\be\label{RZ-eq2} {\real\,}\{e^{i\mu}(1-2ze^{-i\mu}\cos
\nu+z^2e^{-2i\mu})\varphi'(z)\}\geq 0,\;\, z\in \ID. \ee
\end{Thm}

\section{Main Results}

\begin{thm}\label{LiS16-th1}
Suppose that $f_1=h_1+\overline{g_1}\in {\mathcal S}(H^a_{\gamma})$
and $f_2=h_2+\overline{g_2}\in{\mathcal S}(\Omega^b_\beta)$ for some
$a,b\in\ID$.   If $f_1\ast f_2$ is locally univalent in $\ID$, then
$f_1\ast f_2$ is convex in the direction $-\Gamma
:=-(\gamma+\gamma_a+\gamma_b )$, where $\gamma _a=\arg
(1+\overline{a})$ and $\gamma _b=\arg (1+\overline{b})$.
\end{thm}
\bpf First, we recall that $f_1*f_2=h_1*h_2 +\overline{g_1* g_2}$.
In order to prove that $f=h+\overline{g}$ is convex in the direction
of $-\Gamma$, by Theorem \Ref{ThmD}, we need to prove that
$e^{i\Gamma}(h-e^{-2i\Gamma}g)$ is convex in the real direction.

Now, we consider
$$F_1=\left(h_1+e^{-2i(\gamma + \gamma _a)}g_1\right)*\left(h_2-e^{-2i\gamma _b}g_2\right) ~\mbox{ and }~
F_2=\left(h_1-e^{-2i(\gamma +\gamma _a)}g_1\right)*\left(h_2+e^{-2i\gamma _b}g_2\right).
$$
Then we see that
$$\frac{F_1+F_2}{2}=h_1*h_2-e^{-2i\Gamma}g_1* g_2,
$$
where $\Gamma =\gamma+\gamma_a+\gamma_b$, and without loss of
generality, we may assume that $0\leq \Gamma <2\pi$. To prove that
$f_1\ast f_2$ is convex in the direction $-\Gamma$, by Theorem
\Ref{ThmD}, we only need to prove that the function $\Phi$ defined
by
\be\label{LiS16-eq9}
\Phi=e^{i\Gamma} (F_1+F_2)
\ee is convex in the real direction. To do this, we first consider
the following quotients \be\label{LiS16-eq5}
p_1=\frac{h'_1-e^{-2i(\gamma +\gamma _a)}g'_1}{h'_1+e^{-2i(\gamma
+\gamma _a)}g'_1}=\frac{1-e^{-2i(\gamma +\gamma
_a)}\omega_{f_1}}{1+e^{-2i(\gamma +\gamma _a)}\omega_{f_1}} \ee and
\be\label{LiS16-eq5a}
p_2=\frac{h'_2-e^{-2i\gamma _b}g'_2}{h'_2+e^{-2i\gamma _b}g'_2}=\frac{1-e^{-2i\gamma _b}\omega_{f_2}}{1+e^{-2i\gamma _b}\omega_{f_2}},
\ee
where $\omega_{f_j}$ denotes the dilatation of $f_j$ $(j=1,2)$. Clearly,
$p_1$ and $p_2$ are analytic in $\ID$ such that
$$\real\{p_1(z)\}>0 ~\mbox{ and }~ \real\{p_2(z)\}>0~\mbox{in $\ID$.}
$$
We may now set $a'=|1+a|-1$ and $b'=|1+b|-1$. Then, because
$f_1=h_1+\overline{g_1}\in {\mathcal S}(H^a_{\gamma})$, the equation
\eqref{li4-eq2} and the relation \eqref{LiS16-eq5}   give
\be\label{LiS16-eq6} h_1'(z)-e^{-2i(\gamma +\gamma _a)}g_1'(z)
=p_1(z)\left(h_1'(z)+e^{-2i(\gamma +\gamma _a)}g_1'(z)\right)
=\frac{(1+a')p_1(z) }{\left(1-e^{i(\gamma +\gamma _a)}z\right)^2}.
\ee Similarly, since $f_2=h_2+\overline{g_2}\in{\mathcal
S}(\Omega^b_\beta)$, by \eqref{li4-eq3} and
\eqref{LiS16-eq5a}, we have \be\label{LiS16-eq7}
h_2'(z)-e^{-2i\gamma _b}g_2'(z) =p_2(z)\left(h_2'(z)+e^{-2i\gamma
_b}g_2'(z)\right) =\frac{(1+b')p_2(z)}{\left(1+ze^{i(\beta +\gamma
_b)}\right)\left(1+ze^{-i(\beta-\gamma _b)}\right)}. \ee Now,  in
view of the relation
$$zF'_1(z) = \left(h_1(z)+e^{-2i(\gamma +\gamma _a)}g_1(z)\right)*z\left(h_2'(z)-e^{-2i\gamma _b}g_2'(z)\right),
$$
by \eqref{li4-eq2} and \eqref{LiS16-eq7}, we obtain that
\beqq
\frac{F_1'(z)}{(1+a')(1+b')}
&=&\frac{1}{1-e^{i(\gamma +\gamma _a)}z}*\frac{p_2(z)}{\left(1+ze^{i(\beta +\gamma _b)}\right)\left(1+ze^{-i(\beta-\gamma _b)}\right)}\\
&=&\frac{p_2(e^{i(\gamma +\gamma _a)}z)}{\left(1+ze^{i(\beta
+\Gamma)}\right)\left(1+ze^{-i(\beta-\Gamma)}\right)}. \eeqq If we
set $\mu=2\pi-\Gamma$ and $\nu=\pi-\beta$, then easily get
$$(1+ze^{i(\Gamma +\beta)})(1+ze^{-i(\beta-\Gamma)})=1+2ze^{i\Gamma}\cos\beta+z^2e^{2i\Gamma}=e^{i\mu}(1-2ze^{-i\mu}\cos \nu+z^2e^{-2i\mu})e^{i\Gamma}
$$
which by the last equation implies that
\be\label{LiS16-eq8}
{\real}\left \{e^{i\mu}(1-2ze^{-i\mu}\cos
\nu+z^2e^{-2i\mu})e^{i\Gamma}\frac{F_1'(z)}{(1+a')(1+b')}\right\}={\real\,}\{p_2(e^{i(\gamma +\gamma _a)}z)\}> 0
\ee
for $z\in \ID$. Again, because
$$zF'_2(z)=z\left(h_1'(z)-e^{-2i(\gamma +\gamma _a)}g_1'(z)\right)*\left(h_2(z)+e^{-2i\gamma _b}g_2(z)\right),
$$
we have by \eqref{LiS16-eq6},
\beqq
zF'_2(z)
=(1+a')\left(h_2(z)+e^{-2i\gamma _b}g_2(z)\right)*\frac{zp_1(z)}{(1-e^{i(\gamma +\gamma _a)}z)^2}.
 \eeqq
For convenience, we set
$$\phi(z)=h_2(z)+e^{-2i\gamma _b}g_2(z) ~\mbox{ and }~\Psi (z)= \frac{e^{i\Gamma}z}{(1-e^{i(\gamma +\gamma _a)}z)^2}.
$$
Then we observe that
$\phi(z)$ is convex and $\Psi(z)$ is starlike in the unit disk $\ID$. Also, we observe that
$$e^{i\Gamma}\frac{zF'_2(z)}{1+a'} = \phi(z)*p_1(z)\Psi (z)
$$
and
\beqq
\phi(z)*\Psi(z) &=&e^{i\Gamma}z\left(h_2+e^{-2i\gamma _b}g_2\right)'(e^{i(\gamma +\gamma _a)}z)\\
&=&(1+b')\frac{e^{i\Gamma}z}{(1+ze^{i(\beta+\Gamma)})(1+ze^{-i(\beta-\Gamma)})}\\
&=&(1+b')\frac{e^{i\Gamma}z}{1-2ze^{-i\mu}\cos \nu+z^2e^{-2i\mu}}.
\eeqq
Using the last two relations and Theorem \Ref{ThmC}, we find that
\beqq
0<\real\left\{ \frac{\phi(z)*p_1(z)\Psi(z)}{(\phi*\Psi)(z)} \right\} &=& \real\left\{ e^{i\Gamma}z\frac{F'_2(z)}{(1+a')} \cdot \frac {1-2ze^{-i\mu}\cos \nu+z^2e^{-2i\mu}}{(1+b') e^{i\Gamma}z}\right\}\\
&=& \real\left\{e^{i\mu}(1-2ze^{-i\mu}\cos
\nu+z^2e^{-2i\mu})e^{i\Gamma}\frac{F_2'(z)}{(1+a')(1+b')} \right\},
\eeqq which together with \eqref{LiS16-eq8} implies that
$$\real\left\{e^{i\mu}(1-2ze^{-i\mu}\cos
\nu+z^2e^{-2i\mu})e^{i\Gamma}\frac{F_1'(z)+F'_2(z)}{(1+a')(1+b')}\right\}>0~\mbox{
for }~ z\in \ID.
$$
Using Theorem \Ref{ThmE}, we conclude  that $\Phi$ defined by
\eqref{LiS16-eq9}  is convex in the real direction. The proof is
complete.
 \epf

\begin{cor}\label{LiS16-cor1}
Let $a,b\in (-1,1)$. Suppose that $f_1=h_1+\overline{g_1}\in
{\mathcal S}(H^a_{\gamma})$, $f_2=h_2+\overline{g_2}\in{\mathcal
S}(\Omega^b_\beta)$ and $f_1\ast f_2$ is locally univalent in $\ID$.
Then $f_1\ast f_2$ is convex in the direction $-\gamma$.
\end{cor}

It is also possible to generalize Theorems \Ref{ThmA} and \Ref{convolution of
slanted half-plane} to the setting ${\mathcal S}(H^a_{\gamma})$ for
$a\in\ID$.

\begin{thm}\label{LiS16-th2a} 
If $f_k=h_k+\overline{g_k}\in {\mathcal S}(H^{a_k}_{\gamma_k})$ for  $k=1,2$, and $f_1\ast
f_2$ is locally univalent in $\ID$, then $f_1\ast f_2$ is convex in
the direction $-\Gamma$, where
$$\Gamma=(\gamma_1 +\gamma_2+\gamma_{a_1}+\gamma_{a_2}),\;\ \gamma _{a_1}=\arg
(1+\overline{a_1})\;\ \mbox{and}\;\ \gamma _{a_2}=\arg
(1+\overline{a_2}).
$$
\end{thm}
\bpf By Theorem \Ref{ThmD}, it suffices to prove that
$e^{i\Gamma}(h_1*h_2 -{e^{-2i\Gamma}g_1* g_2})$ is convex in the
real direction. We now consider
$$F_1=\left(h_1+e^{-2i(\gamma_1 + \gamma_{a_1})}g_1\right)*\left(h_2-e^{-2i(\gamma_2 + \gamma_{a_2})}g_2\right),
$$
and
$$ F_2=\left(h_1-e^{-2i(\gamma_1 + \gamma_{a_1})}g_1\right)*\left(h_2+e^{-2i(\gamma_2 +
\gamma_{a_2})}g_2\right).
$$
Then we see that
$$\frac{F_1+F_2}{2}=h_1*h_2-e^{-2i\Gamma}g_1* g_2,
$$
which shows that we only need to prove that $e^{i\Gamma}(F_1+F_2)$
is convex in the real direction.

Without loss of generality, we may assume that $0\leq \Gamma <2\pi$. Also, let $\omega_k$ denote the
dilatation of $f_k$  for  $k=1,2$, and consider the following quotients
\be\label{LiS16-eq9a}
q_k=\frac{h'_k-e^{-2i(\gamma_k +\gamma_{a_k})}g'_k}{h'_k+e^{-2i(\gamma_k +\gamma _{a_k}) }g'_k}
=\frac{1-e^{-2i(\gamma_k +\gamma_{a_k})}\omega_{f_k}}{1+e^{-2i(\gamma_k +\gamma
_{a_k})}\omega_{f_k}},\;\ k=1,2,
\ee
which clearly yield that
$$\real\{q_k(z)\}>0 ~\mbox{ for $ z\in\ID$ and for $k=1,2$}.
$$
Because $f_k=h_k+\overline{g_k}\in {\mathcal
S}(H^{a_k}_{\gamma_k})$, \eqref{li4-eq2} and the relation
\eqref{LiS16-eq9a} for $k=1,2$,  give \be\label{LiS16-eq10}
h_k'(z)-e^{-2i(\gamma_k +\gamma _{a_k})}g_k'(z)
=q_k(z)\left(h_k'(z)+e^{-2i(\gamma_k +\gamma _{a_k})}g_k'(z)\right)
=\frac{|1+a_k|q_k(z) }{(1-e^{i(\gamma_k +\gamma _{a_k})}z)^2}. \ee
Since
$$zF'_1(z)=\left(h_1(z)+e^{-2i(\gamma_1 + \gamma_{a_1})}g_1(z)\right)*z\left(h'_2(z)-e^{-2i(\gamma_2 +
\gamma_{a_2})}g'_2(z)\right),
$$
by \eqref{li4-eq2} and the relation \eqref{LiS16-eq10}, we deduce that
$$zF'_1(z)=|1+a_1|\,|1+a_2|\frac{zq_2(e^{i(\gamma_1 +\gamma _{a_1})}z)}{(1-e^{i\Gamma}z)^2}.
$$
Similarly, we can obtain that
$$zF'_2(z)=|1+a_1|\,|1+a_2|\frac{zq_1(e^{i(\gamma_2 +\gamma _{a_2})}z)}{(1-e^{i\Gamma}z)^2}.
$$
Let $\mu=-\Gamma$ and $\nu=0$. Then, we have
$$(1-e^{i\Gamma}z)^2=1-2ze^{i\Gamma}+z^2e^{2i\Gamma}=1-2ze^{-i\mu}+z^2e^{-2i\mu}.
$$
Then the last three equations imply that
$$\real\left\{e^{i\mu}(1-2ze^{-i\mu}\cos
\nu+z^2e^{-2i\mu})e^{i\Gamma}\left(F_1'(z)+F'_2(z)\right)\right\}>0~\mbox{
for }~ z\in \ID.
$$
Using Theorem \Ref{ThmE}, the proof is completed.
 \epf

In the above two theorems, a natural question is to determine conditions on the dilatation of $f_1$ and $f_2$ such that
$f_1\ast f_2$ is locally univalent in $\ID$. Partial answer to this question will be considered in
Section \ref{third-main}.

\section{Convex combination of mappings in ${\mathcal F}^a_{\lambda, \delta}$}\label{sec-main}

Let ${\mathcal F}^a_{\lambda, \delta}$ consist of functions
$F=H+\overline{G}\in{\mathcal H}$ satisfying \be\label{LiS16-eq11}
H'(z)+\overline{\delta^2}e^{-2i\gamma_{a}}G'(z)=\frac{1+a'}{(1+\lambda\delta
e^{i\gamma_a} z)(1+\overline{\lambda}\delta e^{i\gamma_{a}} z)},\ee
where $\lambda,\ \delta\in\partial\ID$, $a\in\ID$ with $a'=|1+a|-1$,
and $\gamma _{a}=\arg (1+\overline{a}).$ If $a$ is real and
$a\in(-1,1)$, then \eqref{LiS16-eq11} reduces to
$$H'(z)+\overline{\delta^2}G'(z)=\frac{1+a}{(1+\lambda\delta
 z)(1+\overline{\lambda}\delta  z)}.
$$
Moreover it is obvious that ${\mathcal S}(H^a_{\gamma})\subseteq {\mathcal
F}^a_{-1, e^{i\gamma}}$ and ${\mathcal S}(\Omega^b_\beta)\subseteq
{\mathcal F}^b_{e^{i\beta}, 1}$.

\begin{lem}\label{convex f after rotation}
Let $F=H+\overline{G}\in{\mathcal F}^a_{\lambda, \delta}$ for
certain values of $\lambda,\ \delta\in\partial\ID$ and $a\in\ID$. Also,  let $\gamma_{a}=\arg (1+\overline{a}).$
Then $F^{\overline{\delta}e^{-i\gamma_a}}=h+\overline{g}$ satisfies that
$$h'(z)+g'(z)=\frac{|1+a|}{(1+\lambda z)(1+\overline{\lambda} z)},
$$
that is, $F^{\overline{\delta}e^{-i\gamma_a}}\in{\mathcal
F}^{a'}_{\lambda, 1}$, where $a'=|1+a|-1$.
\end{lem}
\bpf It follows from the definition of rotation  (see \eqref{LiS16-eq1a}) that
$$h(z)=\delta e^{i\gamma_a} H(\overline{\delta}e^{-i\gamma_a} z)~\mbox{ and }~ g(z)=\overline{\delta}e^{-i\gamma_a} G(\overline{\delta}e^{-i\gamma_a} z),
$$
which yield that
$$h'(z)=H'(\overline{\delta}e^{-i\gamma_a}z) ~\mbox{ and }~\ g'(z)
=\overline{\delta^2}e^{-2i\gamma_a}G'(\overline{\delta}e^{-i\gamma_a}z),
$$
respectively. The last two relations and \eqref{LiS16-eq11} imply that
$$h'(z)+g'(z)=H'(\overline{\delta}e^{-i\gamma_a} z)+\overline{\delta^2}e^{-2i\gamma_a}G'(\overline{\delta}e^{-i\gamma_a}
z)=\frac{|1+a|}{(1+\lambda z)(1+\overline{\lambda} z)}
$$
and the proof is completed. \epf

\begin{thm}\label{convex}
Let $F=H+\overline{G}\in{\mathcal F}^a_{\lambda, \delta}$ for
certain values of $\lambda,\ \delta\in\partial\ID$ and $a\in\ID$. If $F$ is
locally one-to-one and sense-preserving on $\ID$, then $F$ is a
convex harmonic mapping.
\end{thm}
\bpf Let $\gamma _{a}=\arg (1+\overline{a})$. Since
$F^\mu(z)=\overline{\mu}F(\mu z)$ for $\mu\in\partial\ID$, it is
obvious that $F$ is convex if and only if
$F^{\overline{\delta}e^{-i\gamma_a}}$ is convex. Again, since $F$ is
locally one-to-one and sense-preserving on $\ID$, Lemma \ref{the
dilatation after rotation} implies that
$F^{\overline{\delta}e^{-i\gamma_a}}$ is locally one-to-one and
sense-preserving on $\ID$. Let
$F^{\overline{\delta}e^{-i\gamma_a}}=h+\overline{g}.$ By Theorem
\Ref{ThmD}, it suffices to prove that
$\varphi_\theta=:e^{i\theta}h-e^{-i\theta}g$ is convex in the real
direction for all $\theta\in [0, 2\pi)$.

Moreover, Lemma \ref{convex f after rotation} shows that
\beqq\varphi'_\theta(z)&=& e^{i\theta}h'(z)-e^{-i\theta}g'(z)\\
&=& \left(h'(z)+g'(z)\right)\left(i\sin\theta+\frac{h'(z)-g'(z)}{h'(z)+g'(z)}\cos\theta\right)\\
&=&\frac{|1+a|}{(1+\lambda z)(1+\overline{\lambda} z)}
\left(i\sin\theta+\frac{h'(z)-g'(z)}{h'(z)+g'(z)}\cos\theta\right).
\eeqq
Set $\real \lambda=\cos\alpha.$ If
$\theta\in[0,\frac{\pi}{2}]\cup[\frac{3\pi}{2}, 2\pi)$, then we let
$\mu=0$ and $\nu=\pi-\alpha$. This yields that
$$(1+\lambda z)(1+\overline{\lambda}
z)=1+2z\cos\alpha+z^2=e^{i\mu}\left(1-2ze^{-i\mu}\cos\nu+z^2e^{-2i\mu}\right)
$$
and, because $F^{\overline{\delta}e^{-i\gamma_a}}$ is sense-preserving, it follows that
$$\real
\left\{e^{i\mu}\left(1-2ze^{-i\mu}\cos\nu+z^2e^{-2i\mu}\right)\varphi'_\theta(z)\right\}=|1+a|\real
\left\{\frac{h'(z)-g'(z)}{h'(z)+g'(z)}\right\}\cos\theta>0
$$
for $z\in\ID.$
If $\theta\in(\frac{\pi}{2}£¬\frac{3\pi}{2})$, then we let $\mu=\pi$
and $\nu=\alpha$, which yields that
$$-(1+\lambda z)(1+\overline{\lambda}
z)=e^{i\mu}\left(1-2ze^{-i\mu}\cos\nu+z^2e^{-2i\mu}\right)
$$
and thus, we find that
$$\real
\left\{e^{i\mu}\left(1-2ze^{-i\mu}\cos\nu+z^2e^{-2i\mu}\right)\varphi'_\theta(z)\right\}=-|1+a|\real
\left\{\frac{h'(z)-g'(z)}{h'(z)+g'(z)}\right\}\cos\theta>0
$$
for $z\in\ID.$ Using Theorem \Ref{ThmE}, we conclude  that for each
$\theta\in [0, 2\pi)$, the analytic function $\varphi_\theta$ is convex in the real
direction. The proof is complete.
 \epf

In \cite{FHM}, the following result was proved.

\begin{Thm}\label{old convex combination}
{\rm (\cite[Theorem 1.2]{FHM})} Let $\lambda,\
\delta\in\partial\ID$. If $F_j\in{\mathcal F}^0_{\lambda, \delta}$,
$j=1,2,\ldots,n$, then any convex combination of the $F_j$ is a
convex harmonic mapping.
\end{Thm}

In order to prove Theorem \ref{new convex combination} below, as a generalization of Theorem \Ref{old convex
combination}, we need the following lemma.

\begin{Lem}\label{dilation after combination}
{\rm (See \cite{SJW} and also \cite[Lemma 2.1]{FHM})}
Let $\omega_1$ and $\omega_2$ be two analytic functions in the unit
disk that map $\ID$ to itself. Then for any real number $\theta$ and
all $z\in\ID$,
$$\real
\left\{\frac{1-\omega_1(z)\overline{\omega_2(z)}}
{\left(1+e^{-2i\theta}\omega_1(z)\right)\left(1+e^{2i\theta}\overline{\omega_2(z)}\right)}\right\}>0.
$$
\end{Lem}

\begin{thm}\label{new convex combination}
Let $\lambda,\ \delta\in\partial\ID$ and $a\in\ID$. If
$F_j=H_j+\overline{G_j}\in{\mathcal F}^a_{\lambda, \delta}$,
$j=1,2,\ldots,n$, then any convex combination of the $F_j$ is a
convex harmonic mapping.
\end{thm}
\bpf
Let $\omega_j=G_j'/H_j'$ denote the dilatation of $F_j$ for $j=1,2,\ldots,n$, and
$F=H+\overline{G}=\sum^n_{j=1} t_jF_j$, where $t_1, t_2,
\ldots,\ t_n$ are nonnegative real numbers with $\sum^n_{j=1} t_j=1$.
Then
$$H(z)=\sum^n_{j=1} t_jH_j ~\mbox{ and }~G(z)=\sum^n_{j=1} t_jG_j,
$$
which together with \eqref{LiS16-eq11} imply that $F\in{\mathcal
F}^a_{\lambda, \delta}$. Moreover,
$$H'_j(z)=\frac{|1+a|}{(1+\lambda\delta
e^{i\gamma_a} z)(1+\overline{\lambda}\delta e^{i\gamma_{a}}
z)\left(1+\overline{\delta^2}e^{-2i\gamma_{a}}\omega_j(z)\right)},
$$
$$G'_j(z)=\frac{|1+a|\omega_j(z)}{(1+\lambda\delta
e^{i\gamma_a} z)(1+\overline{\lambda}\delta e^{i\gamma_{a}}
z)\left(1+\overline{\delta^2}e^{-2i\gamma_{a}}\omega_j(z)\right)},
$$
$$H'(z)=\frac{|1+a|}{(1+\lambda\delta
e^{i\gamma_a} z)(1+\overline{\lambda}\delta e^{i\gamma_{a}}
z)}\sum^n_{j=1}
\frac{t_j}{\left(1+\overline{\delta^2}e^{-2i\gamma_{a}}\omega_j(z)\right)},
$$
and
$$G'(z)=\frac{|1+a|}{(1+\lambda\delta
e^{i\gamma_a} z)(1+\overline{\lambda}\delta e^{i\gamma_{a}}
z)}\sum^n_{j=1}
\frac{t_j\omega_j(z)}{\left(1+\overline{\delta^2}e^{-2i\gamma_{a}}\omega_j(z)\right)}.
$$
where $\gamma _{a}=\arg
(1+\overline{a})$.
Consider the function
$$\Phi(z)=\left|\sum^n_{j=1} \frac{t_j}{\left(1+\overline{\delta^2}e^{-2i\gamma_{a}}\omega_j(z)\right)}\right|^2-\left|\sum^n_{j=1}
\frac{t_j\omega_j(z)}{\left(1+\overline{\delta^2}e^{-2i\gamma_{a}}\omega_j(z)\right)}\right|^2,\;\
z\in\ID.
$$
Finally, since
$$J_F(z)=\frac{|1+a|^2}{\left|1+\lambda\delta
e^{i\gamma_a} z\right|^2 \left|1+\overline{\lambda}\delta
e^{i\gamma_{a}} z\right|^2}\Phi(z),
$$
Theorem \ref{convex} shows
that we only need to prove $\Phi(z)>0$ in the unit disk. By a
similar reasoning as in the proof of Theorem \Ref{old convex
combination}, we find that $\Phi(z)>0$ in the unit disk. We complete
the proof.
\epf

\section{Further Results on convolution}\label{third-main}

In \cite{LiPo3}, the authors considered certain convolution problems and
proved the following.

\begin{Thm}\label{unit1 in LiPo3}
{\rm (\cite[Theroem 1]{LiPo3})}
Let $f=h+\overline{g}\in {\mathcal S}(H^a_{\gamma})$ with
$$h(z)+e^{-2i\gamma}g(z)=\frac{(1+a)z}{1-e^{i\gamma}z} ~\mbox{ and }~
\omega(z)=e^{2i\gamma}\frac{ze^{i\theta}+a}{1+aze^{i\theta}},
$$
where $\theta\in\IR$ and $a\in (-1,1)$. If one of the following
conditions holds, then $f_0\ast f\in {\mathcal S}_H$ and is convex
in the direction $-\gamma$: \bee
\item [{\rm (1)}] $\cos(\theta-\gamma)=-1$ and $-1/3\leq a<1$.

\item [{\rm (2)}] $-1<\cos(\theta-\gamma)\leq1$ and $a^2<\frac{1}{5-4\cos(\theta-\gamma)}$.
\eee
\end{Thm}

Before we continue further discussion, it is worth to state a remark first
and then recall a couple of results from \cite[Theorems 2 and 3]{LiPo3}.

\br Let $a\in (-1,1)$. For any $f=h+\overline{g}\in\mathcal{H}$, the
representation for $f^a_0$ given by \eqref{LiS16-eq12} quickly gives
that
$$\left(h^a_0\ast h\right)(z)=\frac{(1+a)h(z)+(1-a)zh'(z)}{2},
$$
and
$$\left(g^a_0\ast
g\right)(z)=\frac{(1+a)g(z)-(1-a)zg'(z)}{2}.
$$
Then by a computation, we see that the dilatation
$\widetilde{\omega}$ of $f^a_0\ast f$ is given by
\be\label{li4-eq4}
\widetilde{\omega(z)}=\frac{2ag'(z)-(1-a)z g''(z)}{2h'(z)+(1-a)z
h''(z)}.
\ee
These observations will be helpful in deriving many
convolution theorems.
\er

\begin{Thm}\label{unit2 in LiPo3}\cite[Theorems 2 and 3]{LiPo3}
Let $a\in(-1,1)$ and $f^a_\gamma\in {\mathcal S}(H^a_{\gamma})$ with
the dilatation
$$\omega_{f^a_\gamma}
(z)=e^{2i\gamma}\frac{a-e^{i\gamma}z}{1-ae^{i\gamma}z}.
$$
\bee
\item[{\rm (i)}] Suppose that  $f=h+\overline{g}\in {\mathcal S}(H^0_{\gamma_1})$ with
the dilatation $\omega(z)=e^{i\theta}z^n$, where $n$ is a positive
integer and $\theta\in\IR$. If $a\in[\frac{n-2}{n+2},1)$, then
$f\ast f^a_\gamma$ is convex in the direction $-(\gamma_1+\gamma)$.

\item[{\rm (ii)}] Suppose that $f=h+\overline{g}\in{\mathcal S}(\Omega^0_\beta)$ with the
dilatation $\omega(z)=e^{i\theta}z^n$, where $0<\beta<\pi$,
$\theta\in\IR$ and $n$ is a positive integer. If
$a\in[\frac{n-2}{n+2},1)$, then $f\ast f^a_\gamma$ is convex in the
direction $-\gamma$. \eee
\end{Thm}

Using Theorem \Ref{unit2 in LiPo3}, we have

\begin{thm}\label{unit1}
Let $a_1\in \ID$, $a'_1=|1+a_1|-1$, $0\leq\gamma_1,\;\
\gamma_2<2\pi$, and $f_1=h_1+\overline{g_1}\in {\mathcal
S}(H^{a_1}_{\gamma_1})$ with the dilatation
$$\omega_{f_1}
(z)=e^{2i(\gamma_1+\gamma_{a_1})}\frac{a'_1-ze^{i(\gamma_1+\gamma_{a_1})}}{1-a'_1ze^{i(\gamma_1+\gamma_{a_1})}},
$$
where $\gamma _{a_1}=\arg (1+\overline{a_1})$.

Suppose that $f_2=h_2+\overline{g_2}$ is locally univalent with the dilatation
$$\omega_{f_2} (z)=e^{i \theta}z^n,
$$
where  $n$ is a positive integer and $\theta\in\IR$. Then we have the following:
\begin{enumerate}

\item[{\rm (1)}]
If $f_2\in {\mathcal S}(H^{0}_{\gamma_2})$ and $|a_1+1|\in[\frac{2n}{n+2},2)$, then $f_1\ast f_2$ is convex in the
direction $-(\gamma_1+\gamma_{a_1}+\gamma_2)$.

 \item[{\rm (2)}] If $f_2\in {\mathcal S}(\Omega^0_\beta)$ for some $\beta$ with $0<\beta<\pi$, and
 $|a_1+1|\in[\frac{2n}{n+2},2)$, then $f_1\ast f_2$ is convex in the
direction $-(\gamma_1+\gamma_{a_1})$.
\end{enumerate}
\end{thm}
\bpf   As remarked in the beginning  $f_1$ equals
$f^{a_1}_{\gamma_1}$. The assumption implies that
$a'_1\in[\frac{n-2}{n+2},1)$ and thus, Theorem \Ref{unit2 in
LiPo3}(i) shows that $f_1*f_2$ is locally univalent in $\ID$. The
conclusion follows from Theorem \ref{LiS16-th2a}.

The conclusion of the second part follows from Theorem \Ref{unit2 in
LiPo3}(ii) and a similar reasoning as in the proof of the first
part. \epf

Finally, we state and prove the following result.

\begin{thm}\label{unit3}
For $j=1,2$, let $a_j\in \ID$, $a'_j=|1+a_j|-1$, $\gamma _{a_j}=\arg (1+\overline{a_j})$, $0\leq\gamma_1,\;\
\gamma_2<2\pi$, and $f_j=h_j+\overline{g_j}\in {\mathcal
S}(H^{a_j}_{\gamma_j})$ with the dilatation
$$\omega_{f_1}
(z)=e^{2i(\gamma_1+\gamma_{a_1})}\frac{a'_1-ze^{i(\gamma_1+\gamma_{a_1})}}{1-a'_1ze^{i(\gamma_1+\gamma_{a_1})}}
~\mbox{ and }~
\omega_{f_2} (z)=e^{2i(\gamma_2+\gamma_{a_2})}\frac{a'_2+ze^{i\theta}}{1+a'_2ze^{i\theta}},
$$
If one of the following holds
\bee
\item [{\rm (1)}]   $\cos(\theta-\gamma_2-\gamma_{a_2})=-1$ and
$1+3a'_1+3a'_2+a'_1a'_2\geq0$,

\item [{\rm (2)}] $\cos(\theta-\gamma_2-\gamma_{a_2})=1$ and
$1+3a'_1+3a'_1a'_2+(a'_1)^2a'_2>0$,
\eee
then $f_1\ast f_2$ is convex in the direction $-\Gamma$, where $\Gamma=(\gamma_1
+\gamma_2+\gamma_{a_1}+\gamma_{a_2}).$
\end{thm}
\bpf Set $F_j=H_j+\overline{G_j}=f^{e^{-i(\gamma_j+\gamma_{a_j})}}_j$ for $j=1,2$.
Then Lemma \ref{the slanted half-plane mappings after rotation} implies
that $F_j\in {\mathcal S}(H^{a'_j}_{0})$ with
$$H_j(z)+ G_j(z)=\frac{(1+a'_j)z}{1-z},
$$
and,
$$ \omega_{F_1}
(z)=\frac{a'_1-z}{1-a'_1 z}\;\ \mbox{and}\;\ \omega_{F_2}
(z)=\frac{a'_2+ze^{i(\theta-\gamma_2-\gamma_{a_2})}}{1+a'_2
ze^{i(\theta-\gamma_2-\gamma_{a_2})}}.
$$
The representation given by
\eqref{li4-eq4} implies that the dilatation $\widetilde{\omega}$ of $F_1*F_2$ becomes
$$\widetilde{\omega}=\frac{2a'_1G'_2-(1-a'_1)zG''_2}{2H'_2+(1-a'_1)zH''_2}.
$$
Since
$$G'_2=\omega_{F_2}H'_2 ~\mbox{ and }~
G''_2=\omega'_{F_2}H'_2+\omega_{F_2}H''_2,
$$
the dilatation $\widetilde{\omega}$ takes the form
\be\label{LiS16-eq13}
\widetilde{\omega}=\frac{[2a'_1\omega_{F_2}-(1-a'_1)z\omega'_{F_2}]H'_2-(1-a'_1)\omega_{F_2}
z H''_2}{2H'_2+(1-a'_1)z H''_2}.
\ee
On the other hand, we have
\beqq
H'_2(z)&=&\frac{1+a'_2}{(1-z)^2(1+\omega_{F_2}(z))},\\
H''_2(z)&=&\frac{(1+a'_2)[2(1+\omega_{F_2}(z))-(1-z)\omega'_{F_2}(z)]}{(1-z)^3(1+\omega_{F_2}(z))^2},\\
\frac{zH''_2(z)}{H'_2(z)}&=&\frac{z[2(1+\omega_{F_2}(z))-(1-z)\omega'_{F_2}(z)]}{(1-z)(1+\omega_{F_2}(z))}, ~\mbox{ and }\\
\omega'_{F_2}(z)&=&\frac{(1-a'^2_2)e^{i(\theta-\gamma_2-\gamma_{a_2})}}{\left(1+a'_2
ze^{i(\theta-\gamma_2-\gamma_{a_2})}\right)^2}.
\eeqq
A tedious calculation shows that \beqq
\widetilde{\omega}&=&\frac{[2a'_1\omega_{F_2}-(1-a'_1)z\omega'_{F_2}](1-z)(1+\omega_{F_2})-(1-a'_1)\omega_{F_2}
z[2(1+\omega_{F_2})-(1-z)\omega'_{F_2}]}{2(1-z)(1+\omega_{F_2})+(1-a'_1)z[2(1+\omega_{F_2})-(1-z)\omega'_{F_2}]}\\
&=&\frac{2(a'_1-z)\omega_{F_2}(1+\omega_{F_2})-(1-a'_1)z(1-z)\omega'_{F_2}}{2(1-a'_1
z)(1+\omega_{F_2})-(1-a'_1)z(1-z)\omega'_{F_2}}\\
&=&\frac{2(a'_1-z)(a'_2+ze^{i(\theta-\gamma_2-\gamma_{a_2})})(1+ze^{i(\theta-\gamma_2-\gamma_{a_2})})-(1-a'_1)(1-a'_2)(z-z^2)e^{i(\theta-\gamma_2-\gamma_{a_2})}}
{2(1-a'_1 z)(1+ze^{i(\theta-\gamma_2-\gamma_{a_2})})(1+a'_2
ze^{i(\theta-\gamma_2-\gamma_{a_2})})-(1-a'_1)(1-a'_2)(z-z^2)e^{i(\theta-\gamma_2-\gamma_{a_2})}}\\
&=&-e^{2i(\theta-\gamma_2-\gamma_{a_2})}\frac{t(z)}{t^*(z)}, \eeqq
where
$$t(z)=z^3+c_2z^2+c_1z+c_0 ~\mbox{ and }~\ t^*(z)=1+\overline{c_2}z+\overline{c_1}z^2+\overline{c_0}z^3,
$$
with
\beqq
c_2&=&\frac{a'_1+3a'_2-a'_1a'_2+1}{2}e^{-i(\theta-\gamma_2-\gamma_{a_2})}-a'_1,\\
c_1&=&a'_2e^{-2i(\theta-\gamma_2-\gamma_{a_2})}-\frac{3a'_1+a'_2+a'_1a'_2-1}{2}e^{-i(\theta-\gamma_2-\gamma_{a_2})}, \mbox{ and }\\
c_0&=&-a'_1a'_2e^{-2i(\theta-\gamma_2-\gamma_{a_2})}.
\eeqq

\noindent{{\bf Case 1.} $\cos(\theta-\gamma_2-\gamma_{a_2})=-1$ and
$1+3a'_1+3a'_2+a'_1a'_2\geq0$.}
\smallskip

In this case, since $\cos(\theta-\gamma_2-\gamma_{a_2})=-1$, the constants $c_0,c_1,c_2$ take the form
\beqq
c_2&=&-\frac{3a'_1+3a'_2-a'_1a'_2+1}{2},\\
c_1&=&\frac{3a'_1+3a'_2+a'_1a'_2-1}{2},\\
c_0&=&-a'_1a'_2,
\eeqq
and thus, we have
$$t(z)=(z-1)\left(z^2-\frac{3a'_1+3a'_2-a'_1a'_2-1}{2}z+a'_1a'_2\right).
$$\
Cohn's Lemma (see \cite{RS2002}) shows that $|\widetilde{\omega}(z)|<1$ and the conclusion for this case holds.


\noindent{{\bf Case 2.} $\cos(\theta-\gamma_2-\gamma_{a_2})=1$ and
$1+3a'_1+3a'_1a'_2+(a'_1)^2a'_2>0$.}
\smallskip

In this case, since $\cos(\theta-\gamma_2-\gamma_{a_2})=1$, we find that
\beqq
c_2&=&\frac{-a'_1+3a'_2-a'_1a'_2+1}{2},\\
c_1&=&\frac{-3a'_1+a'_2-a'_1a'_2+1}{2}\\
c_0&=&-a'_1a'_2,
\eeqq
so that
$$t(z)=z^3-\frac{a'_1-3a'_2+a'_1a'_2-1}{2}z^2-\frac{3a'_1-a'_2+a'_1a'_2-1}{2}z-a'_1a'_2.
$$
Now, we let
$$t_1(z)=:\frac{t(z)-c_0t^*(z)}{z}=b_2z^2+b_1z+b_0 ~\mbox{ and }~
t^{*}_1(z)=:\overline{b_2}+\overline{b_1}z+\overline{b_0}z^2.
$$
A calculation yields that
\beqq
b_2&=&1-(a'_1a'_2)^2,\\
b_1&=&\frac{(1-a'_1)\left(1+3a'_2+3a'_1a'_2+a'_1(a'_2)^2\right)}{2},\\
b_0&=&\frac{(1+a'_2)\left(1-3a'_1+3a'_1a'_2-(a'_1)^2a'_2\right)}{2},
\eeqq
and therefore,
\beqq
b^2_2-b^2_0&=&\frac{(1-a'_1)(1-a'_2)}{4}\left(3+a'_2+a'_1a'_2+3a'_1(a'_2)^2\right)\left(1+3a'_1+3a'_1a'_2+(a'_1)^2a'_2\right)\\
&=&\frac{(1-a'_1)(1-a'_2)}{4}(3+a'_2)\left(1+\frac{a'_1a'_2(1+3a'_2)}{3+a'_2}\right)\left(1+3a'_1+3a'_1a'_2+(a'_1)^2a'_2\right),
\eeqq which shows that $|b_2|>|b_0|$ if and only if
$$1+3a'_1+3a'_1a'_2+(a'_1)^2a'_2>0.
$$
Again, we may let
$$t_2(z)=\frac{\overline{b_2}t_1(z)-b_0t^{*}_1(z)}{z}.
$$
Since $b_2$, $b_1$ and $b_0$ are real numbers, the calculation shows
that the zero point $z_0$ of $t_2(z)$ is given by
$$z_0=\frac{-b_1}{b_2+b_0}=:\frac{u}{v}=-\frac{1+3a'_2+3a'_1a'_2+a'_1(a'_2)^2}{3+a'_2+a'_1a'_2+3a'_1(a'_2)^2}.
$$
We obtain that
$$u^2-v^2=-8\left(1-(a'_2)^2\right)\left(1-(a'_1a'_2)^2\right)<0,
$$
which implies that $z_0\in\ID.$ Cohn's Lemma shows that $|\widetilde{\omega}(z)|<1$ and the
conclusion for this case holds.

Finally, the desired conclusion follows from Theorem \ref{LiS16-th2a}. \epf

\subsection*{Acknowledgments}

The work of the first author is supported by NSF of China (No.
11571216), the Construct Program of the Key Discipline in Hunan
Province, the Science and Technology Plan Project of Hunan Province
(No. 2016TP1020) and the Science and Technology Plan Project of
Hengyang City (2017KJ183). The work of the second author is supported  in part by Mathematical Research Impact Centric Support
(MATRICS) grant, File No.: MTR/2017/000367, by the Science and Engineering Research Board (SERB),
Department of Science and Technology (DST), Government of India.



\begin{thebibliography}{99}






\bibitem{Clunie-Small-84} J.~G.~Clunie and T.~Sheil-Small,
Harmonic univalent functions, \textit{Ann. Acad. Sci. Fenn. Ser.
A.I.} {\bf 9} (1984), 3--25.


\bibitem{Do} M.~Dorff,
Convolutions of planar harmonic convex mappings, \textit{Complex
Variables Theory Appl.} {\bf 45} (2001), 263--271.


\bibitem{DN} M.~Dorff, M.~Nowak and M.~Wo{\l}oszkiewicz,
Convolutions of harmonic convex mappings, \textit{Complex Var.
Elliptic Equ. } {\bf 57}(5) (2012), 489--503.

\bibitem{Du} P.~Duren,
Harmonic Mappings in the Plane, Cambridge Tracts in Mathematics,
156. Cambridge Univ. Press, Cambridge, 2004.

\bibitem{FHM} \'{A}. Ferrada-Salas,   R. Hern\'{a}ndez and M. J. Mart\'{\i}n,
On convex combinations of convex harmonic mappings, \textit{Bull.
Aust. Math. Soc.} \textbf{96} (2017), 256--262.

\bibitem{Good02} M.~R.~Goodloe,
Hadamard products of convex harmonic mappings, \textit{Complex
Variables Theory Appl. } {\bf 47}(2) (2002), 81--92.

\bibitem{HengSch70} W. Hengartner and G. Schober,
On schlicht mappings to domains convex in one direction,
\textit{Comment. Math. Helv.} \textbf{45} (1970), 303--314.

\bibitem{HengSch73} W. Hengartner and G. Schober,
A remark on level curves for domains convex in one direction,
\textit{Applicable Analysis} \textbf{3} (1973), 101--106.


\bibitem{HM} R. Hern\'{a}ndez and M. J. Mart\'{\i}n,
Stable geometric properties of analytic and harmonic functions,
\textit{Math. Proc. Phil. Soc.} \textbf{155} (2013), 343--359.

%
%
%



\bibitem{LiPo1} L. Li and S.~Ponnusamy,
Solution to an open problem on convolutions of harmonic mappings,
\textit{Complex Var. Elliptic Equ.} {\bf 58}(12)(2013), 1647-1653.


\bibitem{LiPo2} L. Li and S.~Ponnusamy,
Convolutions of slanted half-plane harmonic mappings,
\textit{Analysis (Munich)} {\bf 33} (2013), 159--176.

\bibitem{LiPo3} L. Li and S.~Ponnusamy, Note on the convolution of harmonic
mappings, \textit{Bull. Aust. Math. Soc.} (2019), 11 pages;
{\tt doi:10.1017/S0004972719000029}

\bibitem{LiuPo} Z. Liu and S.~Ponnusamy,
Univalency of convolutions of univalent harmonic right half-plane
mappings, \textit{Comput. Methods Funct. Theory} {\bf 17}(2) (2017),
289-302.


\bibitem{M} S.~Muir,
Harmonic mappings convex in one or every direction, \textit{Comput.
Methods Funct. Theory} {\bf 12} (2012), 221--239.

\bibitem{samy95} S.~Ponnusamy,
P\'{o}lya-Schoenberg conjecture for Carath\'{e}odory functions,
\textit{J. London Math. Soc.} {\bf 51}(2) (1995), 93--104.





\bibitem{SaRa2013} S. Ponnusamy and A. Rasila, \textit{Planar harmonic and
quasiregular mappings}, Topics in Modern Function Theory (Editors.
St. Ruscheweyh and S. Ponnusamy): Chapter in CMFT, RMS-Lecture Notes
Series No. 19, 2013, pp. 267--333.

\bibitem{PoSi96} S.~Ponnusamy and V.~Singh,
Convolution properties of some classes analytic functions,
\textit{Zapiski Nauchnych Seminarov POMI} {\bf 226} (1996),
138--154; translation in  \textit{J. Math. Sci. (New York)} {\bf
89}(1) (1998), 1008--1020.

\bibitem{RS2002} Q. T. Rahman and G. Schmeisser, Analytic theory of polynomials, London Mathematical
Society Monographs. New Series, Vol. 26, Oxford University Press, Oxford, (2002)

\bibitem{RZ} W. C. Royster and M. Ziegler,
Univalent functions convex in one direction, \textit{Publ. Math.
Debrecen} \textbf{23} (1976), 339--345.

\bibitem{rs1} St.~Ruscheweyh and T.~Sheil-Small,
Hadamard products of schlicht functions and the P\'{o}lya-Schoenberg
conjecture, \textit{Comment. Math. Helv.} {\bf 48} (1973), 119--135.


\bibitem{SJW} Y. Sun, Y. P. Jiang and Z. G. Wang,
On the convex combinations of slanted half-plane harmonic mappings,
\textit{J. Math. Anal.} \textbf{6} (2015), 46--50.
\end{thebibliography}
\end{document}